\documentclass[11pt]{article}

\usepackage[english]{babel}
\usepackage[top=0.6in, bottom=0.7in, left=1.1in, right=1.1in]{geometry}
\usepackage{amsfonts}
\usepackage{amssymb}
\usepackage{relsize}
\usepackage[utf8]{inputenc}
\usepackage{amsmath}
\usepackage{graphicx}
\usepackage[colorinlistoftodos]{todonotes}

\title{On Marshall Hall's Conjecture and Gaps Between Integer Points on Mordell Elliptic Curves}

\author{Ryan D'Mello\\
 Benet Academy, Lisle, IL\\ \\
 {\small Mailing address: Benet Academy, 2200 Maple Avenue, Lisle, IL 60532, USA} \\ {\small e-mail: 5ryan.dmello@benet.org}}

\date{\today}

\begin{document}
\maketitle

\begin{abstract}
For a non-square integer $x\in \mathbb{N}$, let $k_x$ denote the distance between $x^{3}$ and the perfect square closest to $x^{3}$. A conjecture of Marshall Hall states that the ratios $r_x = \sqrt{x}/k_x$, are bounded above. (Elkies has shown that any such bound must exceed 46.6.)
Let \{$x_{n}$\} be the sequence of \textit{"Hall numbers"}: positive non-square integers for which $r_{x_{n}}$ exceeds 1. Extensive computer searches have identified approximately 50 Hall numbers. (It  can be proved that infinitely many exist.) In this paper we study the minimum gap between consecutive Hall numbers. We prove that for all $n$, $x_{n+1} - x_{n} > \frac{1}{5}x_{n}^{\frac{1}{6}}$, with stronger gaps applying when ${x_{n}}$ is close to perfect even or odd squares ($\approx x_{n}^{\frac{1}{3}}$ or $\approx x_{n}^{\frac{1}{4}}$, respectively). This result has obvious implications for the minimum \textit{"horizontal gap"} (and hence straight line and arc distance) between integer points (whose x-coordinates exceed $k^{2}$) on the Mordell elliptic curves $x^{3} - y^{2} = k$ , a question that does not appear to have been addressed.
\end{abstract}

\section{Introduction} 

\textbf{{\small [Appendix A includes a detailed expository coverage of interesting and relevant history and background on Hall's conjecture, its relationship to other math questions/conjectures, topics, and techniques; and  some key applications of these techniques.]}}
\\ \\
For a non-square integer $x\in \mathbb{N}$, let $k_x$ denote the distance between $x^{3}$ and the perfect square closest to $x^{3}$. A conjecture of Marshall Hall states that the ratios $r_x = \sqrt{x}/k_x$, are bounded above. (Elkies [HTTP7] has shown that any such bound must exceed 46.6.) Extensive computer searches and several articles ([Aand], [Calv], [Elks1]) have focused on finding non-square $x\in \mathbb{N}$ for which $r_{x} > 1$, and only about 50 such numbers, which we shall refer to as \textit{Hall numbers}, are known to date (see [HTTP1] for a listing of some of these numbers). The first three such numbers are 2, 5234, and 8158. All members $d$ of the Danilov-Elkies infinite sequence ([Danl], [HTTP6],[Elks1]) have $r_{d}$ slightly greater than 1, and there is no other such infinite sequence known. The highest r-value so far has been found by Elkies (x = 5853886516781223, for which $r_{x} = 46.6$). 
\\ \\
Hence, if Hall's bound does exist, it must be greater than 46.6. Hall had surmised that the bound was 5, based on results obtained from the limited computing power available at the time. The fact that 46.6 is many times larger than 5 has led to skepticism about the tenability of Hall's conjecture. On the other hand, since no number with r-value greater than 46.6 has been found since Elkies' 1999 discovery, it might be premature to entertain serious doubts that the conjecture is true. A weaker form of the conjecture which replaces $\sqrt{x}$ with $x^{0.5 - \varepsilon}$ (for any $\varepsilon$ between 0 and 0.5) is also unresolved. 
\\ \\
No known research appears to have focused on the minimum gap between members of the sequence $\{x_{n}\}$ of Hall numbers. Estimating these gaps would be of interest not only in the context of Hall's conjecture, but also in the context of the "horizontal" (x-value) separation between integer points (whose x-coordinates exceed $k^{2}$) on Mordell elliptic curves $x^{3} - y^{2} = k$, a topic that also does not seem to have been addressed in the literature. It appears that investigations of gaps related to points on Mordell curves have focused on using height or canonical height as the metric for separation.
\\ \\
It would seem very natural and intuitive to ask what pattern the x-gaps of integer points displays. Obviously, knowing the x-gap between two integer points on the curve would also provide estimates of the straight-line (and arc) distance gap between them. 
\\ \\
This paper proves the following theorem: \\ \\
\textbf{Theorem\textit{}}: Let $\{x_n\}$ be the sequence of Hall numbers. Then
\[
  x_{n+1} - x_n > \frac{1}{5}x_n^{\frac{1}{6}}
  \quad\text{for all $n \ge 1$.}
\]
We show that stronger gaps hold when ${x_{n}}$ is close to a perfect even or odd square ($\approx x_{n}^{\frac{1}{3}}$ or $\approx x_{n}^{\frac{1}{4}}$, respectively). Since all Hall numbers below $3\cdot 10^{18}$  are known and meet these gap requirements, we will assume throughout this paper that $x_{n}$ is a very large number (greater than $10^{18}$).
\\ \\
\section{A Few Preliminary Lemmas and Considerations}
In this section, we prove a few simple lemmas and explain the strategy for proving the main result.
\\ \\
\textbf{\textit{Lemma 1}}: Suppose $x^{\frac{3}{2}} = y + f$, where $x$ and $y$ are positive integers and $0 < f < 1$. Then $k = x^{3} - y^{2} = 2fx^{\frac{3}{2}}$ - $f^{2}$. Also,
\[
  r = \dfrac{\sqrt{x}}{k} = \dfrac{1}{2fx  - \dfrac{f^{2}}{\sqrt{x}}}
\] 
\\ 
Proof: Squaring both sides of $x^{\frac{3}{2}} = y + f$ yields $x^{3}$ = $(y + f)^{2}$ = $y^{2} + 2yf + f^{2}$ = $y^{2} + 2f(y + f) - f^{2}$ = $y^{2} + 2fx^{\frac{3}{2}} - f^{2}$. Subtracting $y^{2}$ from both sides then yields the value of k. Calculating r using this value of k and dividing numerator and denominator by $\sqrt{x}$ yields the second part of the lemma.
\\ \\
\textbf{\textit{Lemma 2}}: Suppose $x^{\frac{3}{2}} = y - f$, where $x$ and $y$ are positive integers and $0 < f < 1$. Then $k = y^{2} - x^{3} = 2fx^{\frac{3}{2}}$ + $f^{2}$. Also, 
\[
  r = \dfrac{\sqrt{x}}{k} = \dfrac{1}{2fx  + \dfrac{f^{2}}{\sqrt{x}}}
\] 
\\ 
Proof: The proof is almost identical to that of Lemma 1 (noting that $y^{2} - 2yf + f^{2}$ = $y^{2} - 2f(y - f) - f^{2}$).
\\ \\
\textbf{\textit{Lemma 3}}: If $x$ is a Hall number ($r_{x}$ $>$ 1) and $x^{\frac{3}{2}} = y \pm f$ ($0 < f < 0.5$), then $ 2fx <  1 + \frac{1}{4\sqrt{x}}$ (or f is less than $\approx \frac{1}{2x}$).
\\ \\
Proof: The proof follows directly from the preceding two lemmas, noting that $f^{2}$ will always be less than $\frac{1}{4}$.
\\ \\
Before proving the next set of lemmas, a key step that enables the proof of the main result needs to be introduced. This step requires that $x^{\frac{3}{2}}$ be estimated by first finding the even perfect square closest to $x$. To that end, note that as m ranges over the positive integers, the values $n^{2}$ + a $(1 \leq a \leq 2n)$ and $n^{2}$ - a $(1 \leq a \leq 2n - 2)$, where n = 2m, range over all non-square $x \in \mathbb{N}$. 
\\ \\
The logic is simple: for non-square $x$ $>$ 2, $n^{2}$ is the even square which is closest to $x$. If $x$ is to the right of $n^{2}$, then the first $n^{2}$ + a expression applies, and if to the left of $n^{2}$, then the $n^{2}$ – a expression applies. For example, when m = 2, the intervals (9,16) and (16,25) exhaust all the non-square integers between $3^{2}$ and $5^{2}$. This is depicted in the figure below, where the ranges marked by the dark horizontal lines (excluding the points on the dotted lines) exhaust all non-square $x \in \mathbb{N}$.
\\ 
\begin{center}
FIGURE - 1
\end{center}
\begin{center}
\includegraphics[scale=0.8]{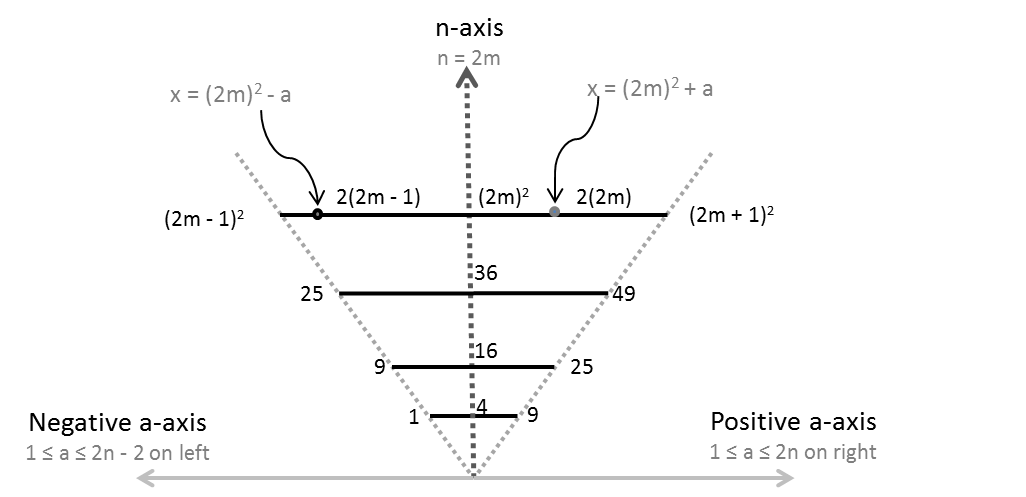} 
\end{center}
The requirement that n be even makes for convenience in some of the calculations that follow. With $x$ expressed as $n^{2} \pm a$ (as above), the Taylor expansion about a = 0 (or, equivalently, the generalized binomial theorem) can be used to expand $x^{3/2}$ = $(n^{2} \pm a)^{\frac{3}{2}}$ = $(n^{2})^{\frac{3}{2}}(1 \pm \frac{a}{n^{2}})^{\frac{3}{2}}$ = $(n^{3})(1 \pm \frac{a}{n^{2}})^{\frac{3}{2}}$. This yields:
\\ \\
$n^{3} + \dfrac{3na}{2} + \dfrac{3a^{2}}{8n} - \dfrac{a^{3}}{16n^{3}} + \dfrac{3a^{4}}{128n^{5}} - \dfrac{3a^{5}}{256n^{7}}$  ......     (for $n^{2} + a$: Case I)
\\ \\
$n^{3} - \dfrac{3na}{2} + \dfrac{3a^{2}}{8n} + \dfrac{a^{3}}{16n^{3}} + \dfrac{3a^{4}}{128n^{5}} + \dfrac{3a^{5}}{256n^{7}}$  ......     (for $n^{2} - a$: Case II)
\\ \\
Let $q_{1}(a)$ = $\dfrac{a^{3}}{16n^{3}} - \dfrac{3a^{4}}{128n^{5}} + \dfrac{3a^{5}}{256n^{7}} - $  ......   (Case I);
\\ \\
and $q_{2}(a)$ = $\dfrac{a^{3}}{16n^{3}} + \dfrac{3a^{4}}{128n^{5}} + \dfrac{3a^{5}}{256n^{7}} + $  ......   (Case II).
\\ \\
It is easy to see that both $q_{1}(a)$ and $q_{2}(a)$ are positive increasing functions of $a$ (their first derivatives are positive, $n$ held constant), the absolute values of the terms in $q_{1}(a)$ and $q_{2}(a)$ are decreasing, the numeric coefficients in the denominators of each term increase by a factor of at least two, and that the difference in the exponents of $a$ and $n$ increases by exactly one from one term to the next. Hence, because $a \leq 2n$ in Case I, $q_{1}(a)$ $\leq$ $q_{1}(2n) < $ 0.5. Similarly, $q_{2}(a)$ $<$ $\dfrac{(a + 1)^{3}}{16n^{3}}$ $<$ 0.5 in Case II (because $(1 \leq a \leq 2n - 2)$.
\\ \\
Suppose $3a^{2}  \equiv  L$ (mod  $8n)$, $0 \leqslant L < 8n$. Then $3a^{2} = 8nt + L$, $t$ being a non-negative integer. Since $n$ is even, $W$ = $n^{3} \pm \dfrac{3na}{2}$ is always an integer. Hence,
\\ \\
$$
x^{\frac{3}{2}} =
\begin{cases}
W + t + \frac{L}{8n} - q_{1}(a)	  & \text{(Case I), } 0 < q_{1}(a)< 0.5 \\
\\
W + t + \frac{L}{8n} + q_{2}(a)	  & \text{(Case II), } 0 < q_{2}(a)< 0.5
\end{cases}
$$
\\
Now, suppose that $x$ is a Hall number. Then, in both cases, $L \neq 0$ because that would imply $3a^{2} \geqslant 8n$ or $ a > \sqrt{n}$, which makes $\frac{1}{16}\frac{a^{3}}{n^{3}}$ (the first term of both $q_{1}(a)$ and $q_{2}(a)$) greater than $1/(16n^{\frac{3}{2}})$, requiring that $q_{1}(a)$ and $q_{2}(a)$ are greater than $1/(32n^{\frac{3}{2}})$  $ > 1/n^{2}$ $ > 1/x$ (for $n > 1,024$), which is impossible by Lemma 3. So, assume $L \neq 0$ in the two paragraphs below.
\\ \\ 
\textbf{For Case 1:} $\frac{L}{8n} < 1$ and, because $x$ is not a perfect square, $\frac{L}{8n}$ cannot equal $q_{1}(a)$. Clearly, if $\frac{L}{8n}$ $>$ $q_{1}(a)$, then $[x^{\frac{3}{2}}]$ is $W + t$, with $\frac{L}{8n} - q_{1}(a)$ being the fractional part of $x^{\frac{3}{2}}$; and if $\frac{L}{8n}$ $<$ $q_{1}(a)$, then  $[x^{\frac{3}{2}}]$ is $W + t - 1$, with $1 + \frac{L}{8n} - q_{1}(a)$ $ = 1 -(q_{1}(a) - \frac{L}{8n})$ being the fractional part of $x^{\frac{3}{2}}$. Thus, from Lemma 3, for $x$ to be a Hall number, $|\frac{L}{8n} - q_{1}(a)|$ would need to be less than $\frac{1 + \frac{1}{4\sqrt{x}}}{2x}$. Note that large values (close to 1) of $\frac{L}{8n}$ would not make $x$ a Hall number because $q_{1}(a) < 0.5$, and so $x^{\frac{3}{2}}$ would remain at least a distance of $\frac{1}{8n}$ $>$ $\frac{1}{x}$ from the nearest integer.
\\ \\
\textbf{For Case 2:} Logic similar to that in the last sentence of the previous paragraph implies that for $x$ to be a Hall number in this case, $\frac{L}{8n} + q_{2}(a)$ would have to be very close to 1, which would require that $|\frac{L}{8n} + q_{2}(a) - 1|$ would have to be less than $\frac{1 + \frac{1}{4\sqrt{x}}}{2x}$ (by Lemma 3). We have thus completed the proof of Lemma 4 below.
\\ \\
\textbf{\textit{Lemma 4}}: If $x$ is a Hall number, then, $L \neq 0$, and in Case I, 
\[
|\frac{L}{8n} - q_{1}(a)| < \frac{1 + \frac{1}{4\sqrt{x}}}{2x}
\]
must hold; whereas in Case II, 
\[
|\frac{L}{8n} - (1 - q_{2}(a))| < \frac{1 + \frac{1}{4\sqrt{x}}}{2x}
\]
must hold.
\\ \\
It is clear from the considerations leading to Lemma 4 that the distance from $x^{\frac{3}{2}}$ to the nearest integer depends on two quantities: $L$ and $q_{1}(a)$ (or $q_{2}(a)$ in Case II), both of which depend on $a$ of course. As $a$ varies from 1 to $2n$ (1 to $2n - 2$ in Case II), $q_{1}(a)$ and $q_{2}(a)$ behave very predictably (being continuous functions of $a$), but $L$ is the wild card which takes on discrete and somewhat unpredictable values, and thus makes the distance between $x^{\frac{3}{2}}$ and the nearest integer difficult to pin down. A very key part of the strategy in the proof is to contain the effects of this unpredictability.
\\ \\
To keep the proofs concise and tidy, we will not prove the lemmas and results separately for each of the two cases (Case I and Case II). We will only address Case I, and point out any noteworthy differences in the proofs for Case II. The logic for both cases is almost identical, with only very minor adjustments needed to adapt Case I proofs to Case II. This would be expected since the Case II numbers are essentially mirror images of Case I numbers ($n^{2}$ being the \textit{mirror}). 
\\ \\
Recall the comment from the Introduction section that $x$ can be assumed to be greater than $3\cdot 10^{18}$, because all the results proved below can be easily verified for the 25 Hall numbers less than $3\cdot 10^{18}$.
\\ \\
\textbf{\textit{Lemma 5}}: If $x = n^{2} + a$ (Case I), $(1 \leq a \leq 2n),$ is a Hall number, then $a > \sqrt[3]{6}(n - 1)^{\frac{2}{3}}$. (For $x = n^{2} - a$ in Case II, the corresponding inequality is $a > 2(n - 1)^{\frac{2}{3}} - 1$.)
\\ \\
Proof: From the equality $\frac{L}{8n} = \frac{L}{8n} - q_{1}(a) + q_{1}(a)$, it follows that $\frac{L}{8n} \leq |\frac{L}{8n} - q_{1}(a)| + |q_{1}(a)|$ or $q_{1}(a)$ = $|q_{1}(a)|$ $\geq  \frac{L}{8n} - |\frac{L}{8n} - q_{1}(a)|$. Using Lemma 4, 
\\ \\
$\dfrac{a^{3}}{16n^{3}}$ $>$ $q_{1}(a)$ $\geq  \dfrac{L}{8n} - |\dfrac{L}{8n} - q_{1}(a)| \geq \dfrac{L}{8n} - \dfrac{1 + \dfrac{1}{4\sqrt{x}}}{2x} > \dfrac{L}{8n} - \dfrac{1 + \dfrac{1}{4n}}{2(n^{2} + a)}$ $> \dfrac{L}{8n} - \dfrac{1 + \dfrac{1}{4n}}{2n^{2}}$
\\ \\
because $n^{2} < x = n^{2} + a$. Multiplying both sides by $16n^{3}$ yields
\\ \\
$a^{3} > 2n^{2}L - \dfrac{16n^{3} + 4n^{2}}{2n^{2}} = 2n^{2}L - (8n + 2)$              
\\ \\
Since $3a^{2}  \equiv  L$ (mod  $8n)$, and because $3a^{2}$ can never be congruent modulo 8 to 1 or 2, the least possible positive value of $L$ is 3. ($L$ cannot be 0 by Lemma 4.) Consequently,
\\ \\
$a^{3} > 2n^{2}L - (8n + 2) \geq 6n^{2}  - 8n - 2 > 6(n - 1)^{2}$, which proves the lemma in Case I.
\\ \\
The proof for Case II parallels that for Case I. The corresponding inequality in this case is 
\\ \\
$\dfrac{(a + 1)^{3}}{16n^{3}} > q_{2}(a) > 1 - \dfrac{L}{8n} - |1 - \dfrac{L}{8n} - q_{2}(a)| > \dfrac{8n - L}{8n} - \dfrac{1 + \dfrac{1}{4\sqrt{x}}}{2x} > \dfrac{8n - L}{8n} - \dfrac{1 + \dfrac{1}{4n}}{2n^{2}}$ 
\\ \\
(invoking Lemma 4). Once again, $L = 0$ would imply that $x$ is not a Hall number. The maximum possible value of $L$ is $8n - 4$ because $3a^{2}$ cannot be congruent modulo 8 to $-1, -2$, or $-3$. So, the minimum possible value of $8n - L$ is 4, and multiplying by $16n^{3}$ yields,
\\ \\
$(a + 1)^{3} > 2n^{2}(8n - L) - (8n + 2) \geq 8n^{2} - 8n - 2 > 8(n - 1)^{2}$,
\\ \\
from which the lemma for Case II follows.
\\ \\
Lemma 5 shows that there are boundaries to the left and right of $n^{2}$ ($n$ even) within which there can be no Hall numbers. A similar result holds in the case of $n^{2}$ ($n$ odd), and this is the focus of Lemma 6.
\\ \\
\textbf{\textit{Lemma 6}}: If $n$ is an odd positive integer and $x = n^{2} \pm a$ is a Hall number, then $a > \sqrt{n}$.
\\ \\
Proof: Note that the previous stipulation for $n$ to be even was solely to ensure that the second term ($\frac{3na}{2}$) in the binomial expansion of $x^{\frac{3}{2}}$ is an integer. Other than that, none of the previous lemmas relied on the evenness of $n$. If $n$ is odd, the integrality of $\frac{3na}{2}$ depends on $a$. If $a$ is an even Hall number, the previous logic works as is and would establish that $a > \sqrt[3]{6}(n - 1)^{\frac{2}{3}}$ in the $n^{2} + a$ case, and that $a > 2(n - 1)^{\frac{2}{3}} - 1$ in the $n^{2} - a$ case. If $a$ is odd, a slight twist to the argument is required although the logic is very similar. 
\\ \\
In the $n^{2} + a$ case ($n$ odd, $a$ odd) the associated expansion (for $x^{\frac{3}{2}}$) is:
\\ \\
$n^{3} + \dfrac{3na}{2} + \dfrac{3a^{2}}{8n} - \dfrac{a^{3}}{16n^{3}} + \dfrac{3a^{4}}{128n^{5}} - ........$ 
\\ \\
The fractional contribution from the second term is 0.5 and from the third term it is $\frac{L}{8n}$. So the only way for $x^{\frac{3}{2}}$ to get close to an integer is for either $\frac{L}{8n}$ to get fairly close to 0.5 or for $\frac{a^{3}}{16n^{3}}$ to get close to 0.5 (with $\frac{L}{8n}$  being small). In the latter case, $a$ would have to be relatively large because $\frac{a^{3}}{16n^{3}}$ being greater that even $ \frac{1}{16} $  would require $a > n$. In the former case, $\frac{L}{8n}$  being greater than even  $\frac{3}{8}$  is possible only when $3a^{2} > 3n$ or $a > \sqrt{n}$. Hence, the $\sqrt{n}$ lower bound for $a$ to be a Hall number would always apply for $n$ odd, regardless of whether $a$ is odd or even.
\\ \\
In the $n^{2} - a$ case ($n$ odd, $a$ odd) the associated expansion (for $x^{\frac{3}{2}}$ ) is:
\\ \\
$n^{3} - \dfrac{3na}{2} + \dfrac{3a^{2}}{8n} + \dfrac{a^{3}}{16n^{3}} + \dfrac{3a^{4}}{128n^{5}} + ........$,
\\ \\ 
and the argument here is similar because the positive terms after the second imply that the fractional part of $x^{\frac{3}{2}}$ can get close to 1 only if either $\frac{L}{8n}$ is greater than $\frac{3}{8}$ or $\frac{a^{3}}{16n^{3}}$ is greater than $ \frac{1}{16} $. The argument in the previous paragraph then applies verbatim, establishing a $\sqrt{n}$ lower bound on $a$ here as well, and thus completing the proof of Lemma 6.
\\ \\
Note that Lemma 6 can be strengthened if one separates odd and even $a$, but we will content ourselves with the $\sqrt{n}$ estimate in the interest of having a more crisp and general articulation of the bound when $n$ is odd. Figure-2 captures Lemmas 5 and 6, and builds on Figure-1, portraying the shaded regions (centered around the perfect squares) free of Hall numbers. 
\\ 
\begin{center}
FIGURE - 2
\end{center}
\begin{center}
\includegraphics[scale=0.75]{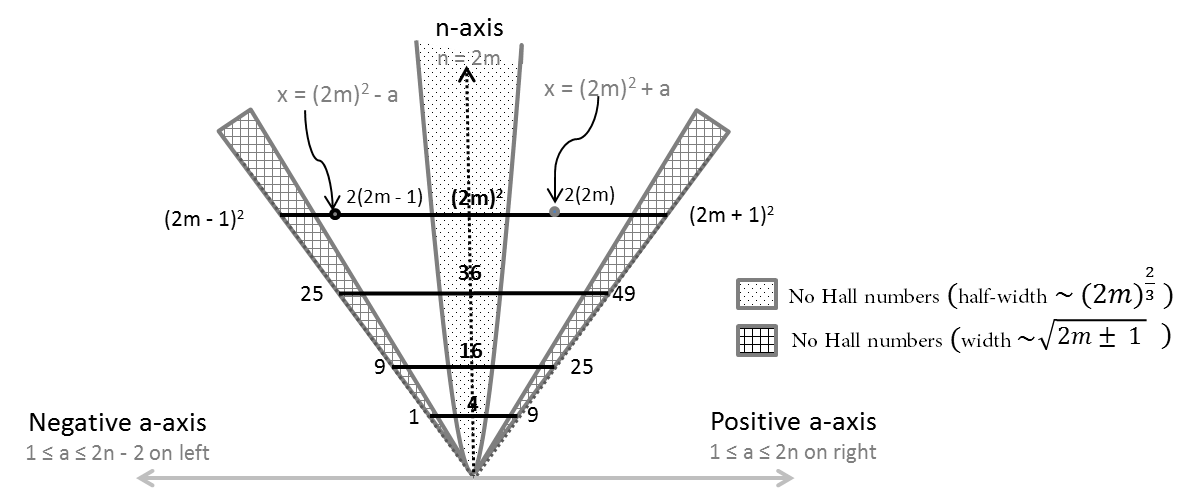} 
\end{center}
\textbf{\textit{Lemma 7}}: If $s$ is a positive integer such that $1 \leqslant a \leqslant  2n - s$ (in Case I), then
\\
\begin{center}
 $\dfrac{3sa^{2}}{16n^{3}}(1 - \dfrac{a}{2n^{2}}) <$  $q_{1}(a + s) - q_{1}(a) <$ $\dfrac{3s(a + s)^{2}}{16n^{3}}$	
\end{center}
The corresponding Case II inequality is
\begin{center}
 $\dfrac{3s(a - s)^{2}}{16n^{3}} <$  $q_{2}(a) - q_{2}(a - s) <$ $\dfrac{3sa^{2}}{16n^{3}}(1 + \dfrac{a}{n^{2}})$
\end{center}
Proof: The proof is a simple application of the mean value theorem. Differentiating both sides of
\begin{center}
$q_{1}(a)$ = $\dfrac{a^{3}}{16n^{3}} - \dfrac{3a^{4}}{128n^{5}} + \dfrac{3a^{5}}{256n^{7}} - $  ......
\end{center}
(with respect to $a$, $n$ held constant) yields
\begin{center}
$q_{1}^{'}(a)$ = $\dfrac{3a^{2}}{16n^{3}} - \dfrac{3a^{3}}{32n^{5}} + \dfrac{15a^{4}}{256n^{7}} - $  ......
\end{center}
The mean value theorem implies that $q_{1}(a + s) - q_{1}(a) = sq_{1}^{'}(a_{0})$, where $a < a_{0} < a + s$. The lemma follows from the fact that $q_{1}^{'}(a)$ is an increasing function of $a$ (because its first derivative is positive), and hence
\\ \\
$\dfrac{3a^{2}}{16n^{3}}(1 - \dfrac{a}{2n^{2}}) = \dfrac{3a^{2}}{16n^{3}} - \dfrac{3a^{3}}{32n^{5}} < q_{1}^{'}(a) < q_{1}^{'}(a_{0}) < q_{1}^{'}(a + s) < \dfrac{3(a + s)^{2}}{16n^{3}}$
\\ \\
For Case II, the proof is almost identical, the corresponding inequality being
\\ \\
$\dfrac{3(a - s)^{2}}{16n^{3}} < q_{2}^{'}(a - s) < q_{2}^{'}(a_{0}) < q_{2}^{'}(a) < \dfrac{3a^{2}}{16n^{3}}(1 + \dfrac{a}{n^{2}})$
\\ \\
One minor but subtle point: Note that for Case II, the $a/\textbf{2}n^2$ from Case I changes to $a/n^2$. The reason for this difference is that when approximating $q_{2}^{'}(a)$ one can easily show that the sum of the third and succeeding terms in the expansion of $q_{2}^{'}(a)$ is much less than the second term.
\\ \\
\textbf{\textit{Lemma 8}}: Suppose $n$ and $s$ are positive integers ($n > 9,709,038$), $s < \frac{1}{4}n^{1/3}$, and suppose $a = kn + f,$ $0 < k < 2.5$, and $|f| < \dfrac{n^{\frac{2}{3}}}{24s}$.  Then (for both Case I and Case II)
\begin{center}
$ \dfrac{3sk^{2}}{2} - \dfrac{4}{5n^{\frac{1}{3}}} <$  $8n|q_{1}(a \pm s) - q_{1}(a)| <$ $ \dfrac{3sk^{2}}{2} + \dfrac{4}{5n^{\frac{1}{3}}}  $
\end{center}
Proof: Substituting $a = kn + f$ in the following inequality from Lemma 7 (for Case I), 
\begin{center}
$\dfrac{3sa^{2}}{2n^{2}}(1 - \dfrac{a}{2n^{2}}) <$  $8n(q_{1}(a + s) - q_{1}(a)) <$ $\dfrac{3s(a + s)^{2}}{2n^{2}}$, we get	
\end{center}
$\dfrac{3sa^{2}}{2n^{2}}(1 - \dfrac{a}{2n^{2}})  = \dfrac{3sk^{2}}{2} + \dfrac{3sfk}{n} + \dfrac{3sf^{2}}{2n^{2}} - \dfrac{3sk^{3}}{4n} - \dfrac{9sk^{2}f}{4n^{2}} - \dfrac{9skf^{2}}{4n^{3}} - \dfrac{3sf^{3}}{4n^{4}}$, and
\\ \\
$\dfrac{3s(a + s)^{2}}{2n^{2}} = \dfrac{3sk^{2}}{2} + \dfrac{3sfk}{n} + \dfrac{3s^{2}k}{n} + \dfrac{3sf^{2}}{2n^{2}} + \dfrac{3s^{2}f}{n^{2}} + \dfrac{3s^{3}}{2n^{2}}$	
\\ \\
Using the bounds for $s$ and $|f|$, and the fact that $3s|f| < \frac{n^{\frac{2}{3}}}{8}$, it follows immediately that the sum of the absolute values of the last five terms in the first expression, and the sum of the last three terms of the second expression are both less than $\frac{4}{n^{\frac{2}{3}}}$, which results in
\\ \\
$ \dfrac{3sk^{2}}{2} + \dfrac{3sfk}{n} - \dfrac{4}{n^{\frac{2}{3}}} <$  $8n(q_{1}(a + s) - q_{1}(a)) <$ $\dfrac{3sk^{2}}{2} + \dfrac{3sfk}{n} + \dfrac{3s^{2}k}{n} + \dfrac{4}{n^{\frac{2}{3}}} $
\\ \\
Note that $\dfrac{3|f|sk}{n} < \dfrac{2.5}{8n^{\frac{1}{3}}}$ and $\dfrac{3s^{2}k}{n} < \dfrac{7.5}{16n^{\frac{1}{3}}}$, and thus the sum of these two terms is less than $\frac{25}{32n^{\frac{1}{3}}}$. The lemma follows from the fact that, for $n > 9,709,038$, $\dfrac{4}{n^{\frac{2}{3}}} < \dfrac{0.6}{32n^{\frac{1}{3}}}$ (and $\frac{25.6}{32} = \frac{4}{5}$). The proof for Case II is identical.
\\ \\ \\
\textbf{\textit{Lemma 9}}: If $L, L^{'}$ and $q, q^{'}$ are any four real numbers, then
\\ \\
$|q - q^{'}| \geq |L - L^{'}| - |L - q| - |L^{'} - q^{'}|$ and $|L - L^{'}| \geq |q - q^{'}| - |L - q| - |L^{'} - q^{'}|$ 
\\ \\
Both inequalities follow from the fact that $L - L^{'} = L - q + q - q^{'} + q^{'} - L^{'}$ and $q - q^{'} = q - L + L - L^{'} + L^{'} - q^{'}$, using the fact that $|a + b + c| \leq |a| + |b| + |c|$ for any real $a$, $b$, and $c$.
\\ \\ 
\textbf{\textit{Lemma 10}}: Let $f(x,y)$ = $k_{0}x^{n}$ + $k_{1}x^{n-1}y$ + .... + $k_{n-1}xy^{n-1}$ + $k_{n}y^{n}$ be a homogeneous polynomial of degree $n \geq 1$ in $x$ and $y$, where the $k_{i}$ are all integers and $k_{n} \neq 0$. Let $a$ and $b$ be non-zero integers such that $k_{n}$ and $a$ are relatively prime and $a^{m}$ divides $f(a,b)$ for some integer $m \geq n$. Then $a \mid b$.
\\ \\
Proof: Let $d$ be the (positive) gcd of $a$ and $b$. Then $a = da_{1}$ and $b = db_{1}$ $(a_{1}$ and $b_{1}$ being relatively prime). Clearly, $f(a,b)$ = $d^{n}f(a_{1},b_{1})$ and
\\ \\
$z$ = $\dfrac{f(a,b)}{a^{m}}$ = $\dfrac{d^{n}f(a_{1},b_{1})}{d^{m}a_{1}^{m}}$ = $\dfrac{f(a_{1},b_{1})}{d^{m-n}a_{1}^{m}}$
\\ \\
Since $z$ is an integer and every term of $f(a_{1},b_{1})$, except possibly the last, is divisible by $a_{1}$, it follows that $a_{1}$ divides $k_{n}b_{1}^{n}$. However, $a_{1}$ is relatively prime to both $k_{n}$ and $b_{1}$, so this forces $|a_{1}|$ = 1, and hence $|a|$ = $d$, which means that $b$ = $db_{1}$ = $|a|b_{1}$, making $b$ a multiple of $a$.
\\ \\
\section{The Main Result}
\textbf{\textit{Theorem 1}}: Suppose $x$ is a Hall number and $s\in \mathbb{N}$, $1 \le s < \frac{1}{5}x^{1/6}$. Then $x \pm s$ cannot be Hall numbers. (The gap $s$ is higher when $x$ is close to a perfect square, per the first paragraph below.)
\\ \\
Proof: From Lemmas 5 and 6 (and Figure 2), it is clear that if $x$ is close to a perfect square (just outside the shaded regions of Figure 2), then the bound for the gap $s$ is actually far better than  $\frac{1}{5}x^{1/6}$, being $n^{\frac{1}{2}} \approx x^{\frac{1}{4}}$ if $x$ is close to an odd perfect square, and $n^{\frac{2}{3}} \approx x^{\frac{1}{3}}$ if $x$ is close to an even perfect square. Hence, we need only consider the cases when $x$ = $n^{2} + a$, $a + s \leq 2n$, $n$ even (for Case I), and $x$ = $n^{2} - a$, $a + s \leq 2n - 2$, $n$ even (for Case II). Also, note that since $x = n^{2} + a < 2n^{2}$, $\frac{x^{\frac{1}{6}}}{5}$ $< \frac{(2n^{2)^{\frac{1}{6}}}}{5}$ $< \frac{n^{\frac{1}{3}}}{4}$ (so the Lemma 8 condition is met). 
\\ \\
The proof for Case II is identical to that for Case I, the only change being that the function $q_{1}$ is replaced by $1 - q_{2}$ (see Lemma 4), and noting that the distance between $1 - q_{2}(a)$ and $1 - q_{2}(a - s)$ is equal to the distance between $q_{2}(a)$ and $q_{2}(a - s)$. For this reason, we will only prove the theorem for Case I. Further, note that proving \textit{"x a Hall number $\Rightarrow$ $x + s$ not a Hall number"} also proves \textit{"x a Hall number $\Rightarrow$ $x -s$ not a Hall number"} (because Lemmas 5 and 6 enable us to assume that $x - s$, $x$, and $x + s$ are all to the right of $n^{2}$). 
\\ \\
As before, set L to be the residue of $3a^{2}$ mod $8n$, so $3a^{2}  \equiv  L$ (mod  $8n)$, $0 \leqslant L < 8n$. Similarly, set $L_{1}$ to the residue of $3(a + s)^{2}$ mod $8n$, so $3(a + s)^{2}  \equiv  L_{1}$ (mod  $8n)$, $0 \leqslant L_{1} < 8n$. This clearly implies that $L_{1} - L$ = $6as + 3s^{2} - 8nv$, for some integer $v$.
\\ \\
We will first prove that if $x + s$ is a Hall number, then $|L_{1} - L|$ must be less than $\sqrt{n}$. This part of the proof is very straightforward. We will then show that a contradiction results if we assume that $|L_{1} - L|$ is less than $\sqrt{n}$, thus resolving the theorem. Proving this contradiction requires more subtle reasoning.
\\ \\
If $x$ and $x + s$ are both Hall numbers, then a straightforward application of Lemma 4 and Lemma 7, and the reasoning in Lemma 9 (noting that $n^{2} < x = n^{2} + a < 2n^{2}$) results in
\\ \\
$|\frac{L}{8n} - \frac{L_{1}}{8n}| \leqslant |\frac{L}{8n} - q_{1}(a)| + |q_{1}(a) - q_{1}(a + s)| + |q_{1}(a + s) - \frac{L_{1}}{8n}|$ 
\\ \\
$< \frac{1 + \frac{1}{4\sqrt{x}}}{2x}$ $+ \frac{1 + \frac{1}{4\sqrt{x+s}}}{2(x+s)} + \frac{3s(a + s)^{2}}{16n^{3}}$ $< \frac{2}{2x} + \frac{2}{2x} + \frac{3s(2n + 2n)^{2}}{16n^{3}}$ $< \frac{2}{x} + \frac{3s}{n}$ $< \frac{2}{n^{2}} + \frac{3x^{\frac{1}{6}}}{5n}$ $< \frac{2}{n^{2}} + \frac{6n^{\frac{1}{3}}}{5n}$  \textbf{...... (A)}
\\ \\
$< \frac{2}{n^{2}} + \frac{2}{n^{\frac{2}{3}}}$ $< \frac{3}{n^{\frac{2}{3}}}$ $< \frac{1}{8\sqrt{n}}$ (for $n > 24^{6} = 191,102,976)$. Multiplying both sides by $8n$ establishes that $|L_{1} - L| < \sqrt{n}$.
\\ \\
We can thus assume going forward that $|L_{1} - L| < \sqrt{n}$. Setting $d = L_{1} - L - 3s^{2}$, and solving the equation $L_{1} - L$ = $6as + 3s^{2} - 8nv$ for $a$ yields
\\ \\
$a = \dfrac{8nv + d}{6s} = n(\dfrac{4v}{3s}) + f$, where $f = \dfrac{d}{6s}$ and $|f| = \dfrac{|d|}{6s} \leqslant \dfrac{|L_{1} - L| + 3s^{2}}{6s}$ $< \frac{\sqrt{n} + \frac{3n^{\frac{2}{3}}}{16}}{6s}$ $< \frac{\frac{4n^{\frac{2}{3}}}{16}}{6s} = \frac{n^{\frac{2}{3}}}{24s}$ (for $n > 16^{6} = 16,777,216$). Note that $v$ cannot be 0 or negative because that would make $a$ less than $\frac{n^{\frac{2}{3}}}{24}$ which contradicts Lemma 5. Also, $k = \frac{4v}{3s}$ must be less than 2.5 or else $a$ would exceed $2n$. Hence, all the conditions of Lemma 8 are met, and it follows that (because $\frac{3sk^{2}}{2} = \frac{8v^{2}}{3s}$)
\\ 
\begin{center}
$ \frac{8v^{2}}{3s} - \frac{4}{5n^{\frac{1}{3}}} <$  $8n(q_{1}(a + s) - q_{1}(a)) <$ $ \frac{8v^{2}}{3s} + \frac{4}{5n^{\frac{1}{3}}}  $  \textbf{  ......... (B)}
\end{center}
There are two possibilities: either $\frac{8v^{2}}{3s}$ is not an integer or it is an integer. We will show that both these possibilities result in a contradiction.
\\ \\
\textbf{Possibility 1:} $\frac{8v^{2}}{3s}$ is NOT an integer. In this case, the denominator of $\frac{8v^{2}}{3s}$ (in reduced form) is at most $3s$, which means that the distance from $\frac{8v^{2}}{3s}$ to the nearest integer is $\geqslant \frac{1}{3s} > \frac{4}{3n^{\frac{1}{3}}}$. Inequality (B) above then requires that $8n(q_{1}(a + s) - q_{1}(a))$ must be at least a distance of $ \frac{4}{3n^{\frac{1}{3}}} - \frac{4}{5n^{\frac{1}{3}}}$ $= \frac{8}{15n^{\frac{1}{3}}}$ $> \frac{1}{2n^{\frac{1}{3}}}$ from the nearest integer. However, Lemma 9 implies, by virtue of inequality (A) above, that (purely on the basis of $x$ and $x + s$ being Hall numbers)  
\\ \\
$|\frac{L}{8n} - \frac{L_{1}}{8n}| - |q_{1}(a) - q_{1}(a + s)| \leqslant |\frac{L}{8n} - q_{1}(a)| + |q_{1}(a + s) - \frac{L_{1}}{8n}| < \dfrac{2}{x}$
\\ \\
Multiplying both sides by $8n$ results in
\\ \\
$|L - L_{1}| - 8n|q_{1}(a) - q_{1}(a + s)| < \frac{16n}{x} < \frac{16n}{n^{2}} = \frac{16}{n}$ $ < \frac{1}{2n^{\frac{1}{3}}}$ (for $n > 182$) \textbf{........ (C)}
\\ \\
which places $8n|q_{1}(a) - q_{1}(a + s)| = 8n(q_{1}(a + s) - q_{1}(a))$ within $\frac{1}{2n^{\frac{1}{3}}}$ of the integer $|L - L_{1}|$, contradicting the earlier conclusion.
\\ \\
\textbf{Possibility 2:} $\frac{8v^{2}}{3s}$ is an integer. It follows from inequalities (B) and (C) above that if $\frac{8v^{2}}{3s}$ is an integer then it must equal $|L - L_{1}|$ (or else $\frac{8v^{2}}{3s}$ and $|L - L_{1}|$ would be separated by at least 1, which is clearly impossible). Substituting $\frac{8v^{2}}{3s}$ for $|L - L_{1}|$ in the equation $a = \dfrac{8nv + d}{6s}$ yields
\\ \\
$a = \dfrac{8nv + d}{6s} = \dfrac{8nv + L_{1} - L - 3s^{2}}{6s} = \dfrac{8nv \pm \dfrac{8v^{2}}{3s} - 3s^{2}}{6s} = \dfrac{24nvs \pm 8v^{2} - 9s^{3}}{18s^{2}}$ 
\\ \\
Since $a$ is an integer and 3 divides the denominator, it follows that $3\mid8v^{2}$ and hence $3\mid v$. Setting $v = 3v_{1}$ and substituting in the above equation yields
\\ \\
$a = \dfrac{72nv_{1}s \pm 72v_{1}^{2} - 9s^{3}}{18s^{2}} = \dfrac{8nv_{1}s \pm 8v_{1}^{2} - s^{3}}{2s^{2}}$ 
\\ \\
This implies that $2\mid s^{3}$, so $s = 2s_{1}$, and another substitution results in
\\ \\
$a = \dfrac{16nv_{1}s_{1} \pm 8v_{1}^{2} - 8s_{1}^{3}}{8s_{1}^{2}} = \dfrac{2nv_{1}s_{1} \pm v_{1}^{2}}{s_{1}^{2}} - s_{1}$
\\ \\ 
Since $a + s_{1}$ is an integer, it follows that $\dfrac{2nv_{1}s_{1} \pm v_{1}^{2}}{s_{1}^{2}}$ is an integer. Using Lemma 10, we can conclude that $s_{1}\mid v_{1}$. Setting $v_{1} = s_{1}v_{2}$ and substituting in the above results in
\\ \\
$a = 2nv_{2} \pm v_{2}^{2} - s_{1} = 2nv_{2} \pm (\frac{v_{1}}{s_{1}})^{2} - s_{1} = 2nv_{2} \pm (\frac{v}{3s_{1}})^{2} - s_{1} = 2nv_{2} \pm (\frac{2v}{3s})^{2} - \dfrac{s}{2}$
\\ \\
Since $k = \frac{4v}{3s}$ is less than 2.5, $(\frac{2v}{3s})^{2}$ is less than 1.6, and so the absolute value of $\pm (\frac{2v}{3s})^{2} - \dfrac{s}{2}$ cannot exceed $1.6 + \dfrac{s}{2} < 1.6 + \frac{n^{\frac{1}{3}}}{8}$. This forces $v_{2}$ to be 1 because $v_{2} \geqslant 2$ would require a to be within a distance of $1.6 + \frac{n^{\frac{1}{3}}}{8}$ from $4n$, which is impossible since $a$ is less than $2n$. (Note that since $v$ is positive, so are $v_{1}$ and $v_{2}$.) However, if $v_{2} = 1$, then the above equation simplifies to 
\\ \\
$a = 2nv_{2} \pm v_{2}^{2} - s_{1} = 2n \pm 1 - s_{1} = 2n \pm 1 - \dfrac{s}{2} - \geqslant 2n - 1 - \dfrac{s}{2}$. Since $s < \frac{n^{\frac{1}{3}}}{4}$, this puts $x = n^{2} + a \geqslant n^{2} + 2n - 1 - \dfrac{s}{2}$ within a distance of $2 + \dfrac{s}{2} < 2 + \frac{n^{\frac{1}{3}}}{8}$ of $(n + 1)^{2}$, which is impossible by Lemma 6 (also see Figure 2) because $n + 1$ is odd.
\\ \\ \\ \\
\textbf{Acknowledgements}
\\ \\
I profoundly thank my parents for a lifetime of love and encouragement. My father inspired in me a deep awe and appreciation for the beauty and applicability of math. He also carefully reviewed this paper, offering many good suggestions to improve readability and clarity, and helped me with the typesetting in LaTeX.
\\ \\
I owe much thanks and appreciation to the late Professor Paul Sally and the University of Chicago, where I attended the Young Scholars Program in the Summer of 2012. Professor Sally planted in my mind the seeds of interest and curiosity about numbers and abstraction. He had an inspiring zest for teaching and an amazing ability to engage his students.
\\ \\
I am very deeply grateful to Professor Joseph Silverman (Brown University) who was the first mathematician to take the time and interest to meaningfully discuss the question I posed to the professional math community (on the MathOverflow website) in August 2014. He shared his  opinions on prior related research, and subsequently offered encouragement and valuable guidance on the publication process. He gave of his time very generously when the paper was being typeset, and offered many very useful LaTeX suggestions for improving format, appearance, and organization. 
\\ \\
Professor Vasudevan Srinivas (Tata Institute of Fundamental Research, Mumbai) provided and sought opinions from his peers regarding the existence of prior research. His frequent communication and suggestions while I was writing up my results were very helpful and motivating. I am very thankful for his time and efforts.
\\ \\
The initial readings that got me curious about the topic of this paper were short articles that I found on the Harvard University web pages of Professor Noam Elkies and Professor Barry Mazur. Professor Noam Elkies very kindly provided detailed answers to questions I asked in an email, informally corroborating my observation that the gap problem addressed here did not appear to be covered in the literature.
\\ \\
Finally, I would like to thank my high school, Benet Academy, for the role it has played in developing my personality, intellect, and interests; and for affording me a diverse set of extracurricular growth opportunities. In particular, I owe much thanks to my teachers and math team coaches in the Benet Math Department for the knowledge they have imparted to me over the years.
\pagebreak
\appendix \begin{center}
{{\Large \textbf{Appendix}\textbf{}}}
\end{center}
\section {History, Background, and Applications}
\subsection{History: The lead-up to Hall's conjecture}
Marshall Hall, Jr. studied math at Yale and Cambridge University, and held faculty positions at several universities. He was a prolific mathematician who made very seminal contributions to group theory, combinatorics, projective geometry, and coding theory. He passed away in London in 1990 en route to a conference to mark his 80{\small th} birthday.
\\ \\
Marshall Hall’s conjecture is concerned with the separation between cubes and squares of integers. If $x$ is an integer (which is not a perfect square), how close can a perfect square get to $x^{3}$? For example, $3^{2}$ = 9 gets to within a distance of one from $2^{3} = 8$. The earliest serious consideration of this separation appears to date back to Euler who, in the eighteenth century, proved that there are no other instances of cubes and squares separated by a distance of one. He used the method of infinite descent, a form of proof by contradiction: assuming a least counterexample exists, and then deducing the existence of an even smaller counterexample.
\\ \\
In 1976, R. Tidjeman ([Tidj]) proved that for positive integer exponents m, n  ($m > 1$ and $n > 1$) there can be only finitely many integer pairs (x,y) such that $x^{m}$ – $y^{n} = 1$, but no such pairs could be identified besides (3,2), with m = 2, n = 3. In 2002, Preda Mihailescu ([Mihl)] finally settled the question definitively by proving that with the exception of (3,2), there are no other (x,y) pairs for which $x^{m}$ – $y^{n}$ $= 1$. 
\\ \\
If, as Euler proved, there are no other cubes besides 8 that get to within a distance of one from a perfect square, how close can a cube (which is not a perfect square) get to a perfect square? In other words, if
\\ \\
$x^{3}$ – $y^{2}$ = k
\\ \\
where x and y are integers, and x is not a perfect square, then can one estimate a lower bound for $|k|$ in terms of some function of x? (Euler’s result implies that the lower bound has to be at least 2, for x $> 2$). The earliest work on this question can be traced to a letter that S. Chowla sent to B. J. Birch (dated 29 September, 1961, as mentioned in [HTTP7]). Chowla investigated the polynomial family 
\\ \\
$x(t) = \dfrac{t(t^{9} + 6t^{6} + 15t^{3} + 12)}{9}$
\\ \\
$y(t) = \dfrac{t^{15}}{27} + \dfrac{t^{12} + 4t^{9} + 8t^{6}}{3} + \dfrac{5t^{3} + 1}{12}$
\\ \\
$k(t) = x^{3}(t) - y^{2}(t) = -\dfrac{3t^{6} + 14t^{3} + 27}{108}$
\\ \\
(where $t$ $\equiv$ $3$ modulo 6, so that x(t) and y(t) are integers). Since x is of the order of $t^{10}$ and k(t) is of the order of  $t^{6}$, the above identity establishes that there are infinitely many (x,y) pairs for which $|k|$ $< Cx^{\frac{3}{5}}$, C being a constant. 
\\ \\
The obvious question would then be whether a similar result could be proved with an exponent of x less than $\frac{3}{5}$. Hall’s experimental evidence for values of x up to $700,000$ seemed to indicate that k is always greater than $C\sqrt{x}$ for some constant C.  For x = 5,234 (hence $\sqrt{x}$ $\approx$ 72.35) and y = 378,661, $x^{3}$ – $y^{2}$ = -17. 
Since 17/72.35  $\approx$ 0.235, it follows that for Hall’s observation to be generalized C would have to be less than 0.235 or, equivalently, r = $\dfrac{\sqrt{x}}{k}$  would have to be greater than 4.26.
\\ \\  
Davenport ([Davn]) proved in 1965 that for polynomials f(t) and g(t) with complex coefficients, and such that $f(t)^{3}$ $\neq$ $g(t)^{2}$, the degree of $f(t)^{3}$ – $g(t)^{2}$ is always greater than or equal to $\dfrac{degree(f(t))}{2}$ + 1. Marshall Hall’s conjecture was prompted by these results and stated:
\\ \\
\textbf{Marshall Hall's Conjecture\textit{}}: If $x^{3} - y^{2} = k$ for integers x, y (x not a perfect square), then $|k| > \dfrac{\sqrt{|x|}}{5}$ (equivalently, $r = \dfrac{\sqrt{|x|}}{|k|} < 5)$
\\ \\
The measure r indicates how small $|k|$ is relative to $\sqrt{|x|}$ : the higher the value of r, the smaller the value of $|k|$ relative to $\sqrt{|x|}$ . 
\subsection{Subsequent efforts}
Hall’s conjecture sparked a search for integer pairs (x,y) for which $|k|$ = $|x^{3}$ – $y^{2}|$ is less than $\sqrt{|x|}$ (or $r > 1$), and $x$ is not a perfect square. For convenience, we shall refer to such a number $x$ as a \textit{Hall number}. As mentioned in [Elks1], the Danilov-Elkies Fermat-Pell family
\\ \\
	$x$ = $5^{5}t^{2} - 3000t + 719$
\\
	$Y$ = $(5^{3}t^{2} – 114t + 26)(5^{6}t^{2} – 5^{3}(123)t + 3781)^{2}$
\\
    $k$ = $27(2t -1)$
\\ \\
satisfies the identity $x^{3}$ – $Y$ = $k$
\\ \\
It can be proved using the theory underlying solutions to Pell’s equation that $5^{3}t^{2} – 114t + 26$ is a perfect square for infinitely many values of $t$. ($t$ = -5 and $t$ = -10,150,883 are the smallest such $t$, in absolute value, and the corresponding values of $x$ are 93844  and 322001299796379844.) For these values of $t$, $Y$ would be a perfect square (say $y^{2}$). Hence, there are infinitely many (x,y) pairs for which $x^{3}$ – $y^{2}$ = $k$ and $|k|$ $\approx C\sqrt{x}$, C = $\dfrac{54}{5^{\frac{5}{2}}}$ $\approx$ 0.966. Equivalently, $r = \sqrt{|x|}/|k| > 1/C$ ~ $\approx$ 1.035. The five values of $x$ calculated to date are documented in [HPPT6]. They are the two mentioned above (93844  and 322001299796379844) and:
\\ \\
1114592308630995805123571151844,  3858108676488182444301031186675778188809844, and
\\
13354661111806898918013326915229994453818137920195953844.
\\ \\
None of the pairs in the Danilov-Elkies family would violate Hall’s conjecture. To disprove the conjecture, a value of $r$ above 5 would be required.
\\ \\
In 1998, Elkies ([Elks1]) developed an algorithm to search for values of $x < X$ with $r$ $>$ 1. He encoded it in a C program using 64-bit integer arithmetic that ran in $\sqrt{X}log(X)$ time. His program found all values of $x$ less than $X = 3(10^{18})$. Two values had $r$ exceeding 5. Specifically, 38115991067861271 had $r$-value 6.5 and 5853886516781223 had $r$-value 46.6. This disproved Hall’s conjecture in its C = 1/5 (or $r$ = 5) formulation. In the same paper, Elkies proved that there can be no more than $\sqrt{X}log(X)$ Hall numbers below $X$.
\\ \\
Subsequent searches focused on finding larger Hall numbers. The interested reader can refer to [Aand], [Calv], and [Elks1]. [Aand] lists all known Hall numbers, 25 of which are greater than $10^{20}$ (and two among these are greater than $10^{40}$). Note, however, that the Hall numbers above 3($10^{18}$) were identified using a probabilistic search algorithm which was not guaranteed to find all Hall numbers in a range, but only those that met specific criteria (which increased the likelihood of their being Hall numbers). Hence, there likely are many Hall numbers in the $10^{20}$ to $10^{40}$ range that have not been identified. The 46.6 value for $r$ found by Elkies still stands as a record.
\\ \\
While Hall’s conjecture in its  $r = 5$ formulation has been disproved, it could still hold with a value of $r$ greater than 46.6 (or a value of C less than 1/46.6 = 0.0215). Hence, Marshall Hall’s conjecture remains unresolved to this day.
\\ \\
A weaker form of Hall’s conjecture, formulated by Stark and Trotter around 1980, states that for any $\varepsilon > 0$, there is a constant C($\varepsilon$), depending only on $\varepsilon$, such that for integers $x, y$ (where $y^{2} \neq  x^{3}$)
\\ \\
$|k| = |x^{3}  -  y^{2}| > C(\varepsilon)x^{\frac{1}{2} - \varepsilon}$ whenever $x >$ C($\varepsilon$)
\\ \\ 
The fact that even this weaker version is unresolved indicates that there is a fundamental gap in our current understanding of the separation between cubes and squares of integers. A generalization of Hall's conjecture is Pillai's conjecture which is associated with separation of other distinct powers of integers (like fourth and fifth powers). It is also unresolved.
\\ \\
\subsection{Hall's Conjecture and Elliptic Curves}
Hall’s conjecture originated from Mordell’s equation $y^{2} = x^{3} + k$ associated with an elliptic curve sometimes known as Mordell’s curve. A more general form of an elliptic curve is the set of points (x,y) in the plane that satisfy an equation of the form: $y^{2} = x^{3} + ax + b$, ($a$ and $b$ are constants.) The \textit{"plane"} here is defined by an underlying field (for example, the complex, real or rational numbers), and $a$ and $b$ are in this field. By including a point at infinity, an addition operator can be defined on the points on the curve. 
\\ \\
To add two points $P_{1}$ and $P_{2}$, one would extend the line joining them until it intersects the curve at a third point, and then set the sum of $P_{1}$ and $P_{2}$ to the reflection (in the x-axis) of this third point. (To add a point to itself, the tangent line would be used instead.) This addition is obviously commutative, but it takes a bit of computation to show that it is also associative. There is an additive identity (which is the point at infinity), and each point has an additive inverse (its reflection in the x-axis). Hence, the points form a group under this addition, an observation that is credited to Henri Poincare. 
\\ \\
If the underlying field contains the rational numbers, it is easy to prove that the points (x,y) on this curve for which x and y are both rational form a subgroup of this group (assuming that $a$ and $b$ are rational). In 1922, Louis Mordell proved that this subgroup is finitely generated ([Mord]). (In his 1928 doctoral dissertation, Andre Weil ([Weil]) generalized this result to include certain other types of fields - so this result is referred to as the Mordell-Weil theorem.) Incidentally, the proof uses a form of infinite descent (which, as was pointed out at the beginning of this appendix, was used by Euler when investigating cubes and squares that are separated by a distance of one).
\\ \\
The fact that this subgroup is finitely generated implies that it is a sum of a finite number of finite cyclic groups and a finite number of copies of $\mathbb{Z}$. The subgroup of this group comprising the finite sum of finite cyclic groups is called the torsion subgroup, and elements of the torsion subgroup are referred to as torsion points. Clearly, they are the points P of finite order: meaning that by computing P, $P + P$, $P + P + P$, etc., we will sooner or later reach P again. In 1977, Barry Mazur ([Mazr1]) proved that the torsion subgroup must be one of only 15 possible groups, the largest of which has 16 elements. 
\\ \\
If an elliptic curve is defined by an equation with integer coefficients, then the Nagell-Lutz theorem states that any torsion point which is rational must actually be an integer point; and in 1928, Carl Siegel proved that an elliptic curve can have only finitely many integer points. (Note that an integer point may not be a torsion point.)
\\ \\
The number of copies of $\mathbb{Z}$ in the decomposition of the subgroup of rational points is referred to as the rank of the elliptic curve. Most curves have small rank, and the elliptic curve with the highest exactly known rank (of 19) is one discovered by Noam Elkies, who also discovered a curve with rank at least 28. (Some examples of curves having high rank are at [HTTP11].) It is not known whether or not there are curves with rank 29 or higher. In 1987, Joseph Silverman ([Silv]) proved that the number of integer points on the Mordell curve C: $y^{2} = x^{3} + k$ (k being $6^{th}$ power free) cannot exceed $K^{1 + rank(C)}$ for some absolute constant K which is independent of C. (Rank(C) is the rank of the subgroup of rational points.)
\\ \\
Determining which numbers can be possible ranks of an elliptic curve is a very difficult open problem. The Birch and Swinnerton-Dyer conjecture, one of the six unresolved Millennium problems (see [HTTP12]), each carrying a million dollar prize, is related to this open problem. 
\\ \\
Interestingly, an elliptic curve called the Frey elliptic curve played a key role in  the proof of \textit{Fermat's Last Theorem} which states that the equation $x^{n} + y^{n} = z^{n}$ $(n \geqslant 3, n \in \mathbb{N})$ cannot have any integer solutions. In 1982, Gerhard Frey demonstrated that if Fermat's Last Theorem had a counterexample, then the Frey elliptic curve would not have certain properties. Andrew Wiles spent six years proving a result which implied that the Frey curve would have the properties in question, and hence ruling out the possibility of a counterexample to Fermat's Last Theorem. 
\\ \\
A generalization of Fermat's Last Theorem, known as Beal's Conjecture (and also as the Tijdeman-Zagier conjecture), states that the equation $x^{l} + y^{m} = z^{n}$ $(l, m, n \geqslant 3;$ $l, m, n \in \mathbb{N})$ has no integer solutions with $x, y, z$ relatively prime. (Infinitely many solutions exist if the \textit{relatively prime} condition is dropped, example $3^{3} + 6^{3} = 3^{5}$.) Andrew Beal, a billionaire banker and amateur mathematician, is offering a million dollars for a proof or counterexample (see [HTTP8]). 
\\ \\
The Beal conjecture has been proved for special cases of $(l,m,n)$ using methods that extend the theory of elliptic curves to a broader class of curves. In the context of this conjecture and the broader class of curves (whose \textit{genus} is greater than one) two theorems stand out: Falting's theorem and the Darmon-Granville theorem. The latter theorem implies that for any specific $(l,m,n)$, there can be at most finitely many relatively prime $x, y, z$ solutions. (Refer to [Elks2] and [Mazr2] for interesting coverage of these topics.)
\\ \\
When a = 0, the elliptic curve represents a Mordell curve. A few examples of Mordell curves are shown below (extracted from [HTTP9]). The black dots on the curve indicate lattice points (both $x$ and $y$ are integers). Sloan’s sequence ([HTTP10]) lists some values of $k$ (6, 7, 11, 13, 14, 20, 21, 23, 29, 32, 34, 39, 42) for which the Mordell curve $C: y^{2} = x^{3} + k$ has no lattice points.
\\ \\
\includegraphics[scale=1]{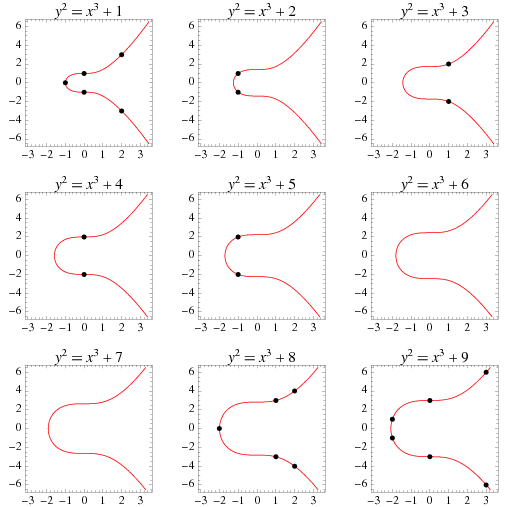}  
\\ \\
Hall’s conjecture has implications for the integer solutions to Mordell’s equation: $y^{2} = x^{3} + k$. As noted earlier, Siegel’s theorem states that this equation (for $k \neq 0$) can have at most finitely many solutions. Baker  proved that $max(|x|,|y|)$ $<$ $e^{c|k|^{1000}}$, where c = $10^{10}$. Stark improved on this result and was able to show that for any $\varepsilon > 0$ there is a constant $c = c(\varepsilon)$ (depending only on $\varepsilon$) such that $max(|x|,|y|)$ $<$ $e^{c|k|^{1 + \varepsilon}}$. However, both of these estimates are exponential. 
\\ \\
Suppose Hall’s conjecture were true, say for $r = 50$. This would imply that $\frac{\sqrt{|x|}}{|k|}$ $< 50$ or $\sqrt{x} < 50k$ or $x < 2500k^{2}$.  Consequently, $y^{2} = x^{3} + k$ $< (2500k^{2})^{3}$ + $k =$ $50^{6}k^{6} + k$ $< 50^{6}(|k| + 1)^{6}$. So, $y < 50^{3}(|k| + 1)^{3}$. Hence, Hall’s conjecture, if true, would imply that $max(|x|,|y|)$ is bounded by a polynomial in $k$, which would be a far better bound than that of Stark.
\\ \\
\subsection{Applications of Elliptic Curves}
Elliptic curves have applications in public-key cryptography. Public-key cryptography is used to encode data in a manner that is hard to decipher, and it is thus applicable in situations where sensitive data needs to be stored and/or transmitted. Examples are defense/intelligence agency communication, password information and credit card transactions on the internet, data storage in the cloud, intra-company communication of proprietary information, transmission of private information (like health care records), etc.  
\\ \\
Public-key cryptography requires two keys: a private key and a public key. The public key is used for encryption and the private key is used for decryption. While the public key and the private key have a mathematical link, it is essential to ensure that determining the private key from the public key is virtually impossible. The goal is to make the deciphering of the private key from the public key not computationally feasible. In other words, there would just be so many possibilities to check, that with today’s computing speeds it would take several decades of computing time to make this determination possible.
\\ \\
The first public-key cryptographic system, which is still being widely used, is known as RSA (named after the initials of the three developers of the underlying algorithm: Ron Rivest, Adi Shamir and Leonard Adleman). The difficulty in deciphering the private key essentially stems from the fact that it is easy to get a computer to multiply two large primes, but difficult to have it factor the product. For example, suppose one wanted to run a simple divisibility check on all numbers between $10^{20}$ and $10^{30}$, and suppose each check can be done in a billionth of a second. This would take several billion centuries, by which time the sun would have burned out!
\\ \\
Elliptic curve cryptography (ECC) (which was introduced in 1985) has two advantages: 
\\ \\
1. It has a smaller public key size which reduces storage and bandwidth/transmission time requirements. For example, a 256-bit ECC public key provides security which is comparable to a 3072-bit RSA public key; and 
\\ \\
2. Deciphering the ECC public key is much harder  (more computationally intensive) than deciphering an RSA public key. The following excerpt from [HTTP3] puts this comparison in  perspective:
\\ \\
		\begin{scriptsize}
		\textit{To visualize how much harder it is to break, Lenstra recently introduced the concept of "Global Security." You can compute how much energy is needed to break a cryptographic algorithm, and compare that with how much water that energy could boil. This is a kind of cryptographic carbon footprint. By this measure, breaking a 228-bit RSA key requires less energy than it takes to boil a teaspoon of water. Comparatively, breaking a 228-bit elliptic curve key requires enough energy to boil all the water on earth. For this level of security with RSA, you'd need a key with 2,380-bits.
		}
		\end{scriptsize}
\\ \\
More information on the documented advantages of ECC over RSA can be found at [HTTP3], and at the US National Security Agency website: The Case for Elliptic Curve Cryptography ([HTTP2]).
\\ \\
ECC uses elliptic curves over a finite field. The simplest example of a finite field is the set of all integers modulo p (under addition and multiplication), where p is a prime. The elliptic curve is then defined as the set of all points (x,y), where x and y are in the field and satisfy
$y^{2} = x^{3} + ax + b$
\\ \\
(which, when a = 0, is just the Mordell curve $y^{2} = x^{3} + b$).
\\ \\
For example, in the field of integers modulo 7 (under addition and multiplication) the point (2,2) satisfies the equation
\\ \\
$y^{2} = x^{3} + x + 1$
\\ \\
because $2^{2}$ $\equiv$ 4 (mod 7) and $2^{3} + 2 + 1 = 8 + 2 + 1 =  11$ $\equiv$ 4 (mod 7). 
\\ \\
The set of all such points (x,y) on the curve (along with the point at infinity) forms a commutative group, with the point at infinity acting as the identity element of this group. It is easy to show that a finite field has $p^{n}$ elements for some prime $p$ (and positive integer $n$), so the overall space in which an elliptic curve over a finite field is defined will have $p^{2n} + 1$ points for some prime p (the 1 being added to include the point at infinity). Some illustrative examples can be found at [HTTP4].
\\ \\
Because it is so difficult to decipher, the adoption of ECC is growing. Governments and defense organizations, including the United States government, are utilizing this relatively newer type of cryptography, making it nearly impossible to tamper with internal communications. Companies are beginning to apply elliptic curve cryptography as well. For example, Apple implements ECC in order to provide signatures in Apple’s iMessage service. Several other text messaging or Short Message Service (SMS) apps incorporate ECC. For example, searching for ECC on the Android App Market (now called Google Play) will display several examples (or just Google a term like \textit{Google Play elliptic curve cryptography}). 
\\ \\
The increasing use of remote patient monitoring systems via the use of embedded mobile devices has raised concerns about privacy in the transmission of patient information from such systems. Once again, ECC has been recommended as an efficient and effective solution by several health care and security experts ([Malh]). 
\\ \\
Similarly, wireless enterprise-wide instant messaging (IM) facilitates collaboration and good communication within a company, but there is always the concern that these messages could be intercepted by wily competitors who could steal valuable proprietary information. Once again, ECC has been deployed to address these concerns.
\\ \\
It is also being used as a measure to secure cloud computing. ECC is transforming secure data storage and communications, and is turning into a trustable solution for online security ([Gupt]). Bitcoin uses an ECC-based digital signature algorithm to ensure that only rightful owners can access and spend their funds. Internet hosting companies are incorporating ECC into their Secure Socket Layer (SSL) technology.
\\ \\
Some of the more recent versions of well-known web browsers, such as Google Chrome and Firefox, also incorporate ECC-based algorithms for communicating securely. If you are using a recent version of Google’s Chrome browser, typing in the URL for any secure site (one beginning with https) and clicking the lock symbol (to the left of the URL) will show the encryption scheme being used. For example, go to https://blog.cloudfare.com and click on the lock symbol. You will see something like the following:
\\
\begin{center}
\includegraphics[scale=0.52]{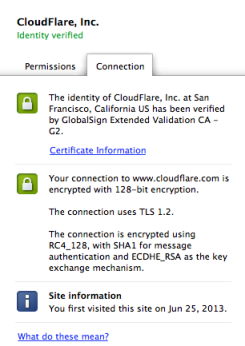} 
\end{center}
Here, the ECDHE stands for Elliptic Curve Diffie Hellman Ephemeral which incorporates an algorithm based on elliptic curves. (The RSA component means that RSA-based algorithms are used to establish the identity of the server.)
\\ \\
\subsection{Hall’s conjecture and the ABC conjecture}
The abc conjecture (also known as the Oesterlé–Masser conjecture) was proposed by Joseph  Oesterlé and David Masser in the 80’s. It is one of the most well-known and challenging open questions in math and its resolution would settle  many outstanding unsolved math problems including the weak form of Marshall Hall’s conjecture. 
\\ \\
The key concept underlying the abc conjecture is simple. Let us define a triple to be an ordered triad (a,b,c), where a and b are co-prime positive integers and a +  b = c. We define the radical r(a,b,c), or just r for short, as the product of the \textbf{distinct} primes dividing abc. Consider the triple (125,3,128).
\\ \\
$125 + 3 = 128$ or $5^{3} + 3 = 2^{7}$
\\ \\
In this case r is 5 x 3 x 2 = 30, which is less than $c = 128$. Triples which have this property $(r < c)$ are referred to as abc triples. Experimenting with a few examples will quickly show that abc triples are quite rare (a randomly chosen triple is very likely to not be an abc triple). For example, there are only 15 abc triples with c less than 300. It is very easy to prove that there are infinitely many abc triples: the family ($1, 9^{n} – 1, 9^{n}$), where n is a positive integer, contains only abc triples because $9^{n} – 1$ is always divisible by 8, hence r would always be less than 0.75($9^{n} – 1$).
\\ \\
How small can the radical r of an abc triple get relative to c? One measure would be to calculate the ratio c/r. The greater this ratio, the smaller r is in comparison to c. However, a more convenient measure is to calculate the ratio of the logarithms: $q =  \frac{log(c)}{log(r)}$. (This ratio is independent of the base on account of the change of base formula.) We refer to q as the quality of the abc triple.
\\ \\
The highest known value of q is about 1.62991 associated with the triple $(2,109.3^{10},23^{5})$, which is known as Reysatt’s triple. The second highest value of q is about 1.62599 and is associated with the triple $(11^{2},3^{2}5^{6}7^{3},23.2^{21})$, known as Benne de Weger’s triple. As of 2013, there were only 140 known triples whose quality was at least 1.4.
\\ \\
The abc conjecture states that for any $\varepsilon >$ 0 there are only finitely many abc triples whose quality q exceeds 1 + $\varepsilon$. Given that we defined q as the ratio of logs, this statement is obviously equivalent to the statement that for any $\varepsilon$ $> 0$ there are only finitely many abc triples such that
$c > r^{1 + \varepsilon}$. It can be proved that the abc conjecture implies the weak form of Hall's conjecture. It does not imply the strong form of Hall's conjecture. [Wald] has some interesting and useful coverage on the abc conjecture.
\\
\subsection{Wrap-up}
As demonstrated in the preceding appendix subsections, Hall's conjecture is not merely a stand alone curiosity, but exists in a larger ecosystem of interesting problems, techniques, conjectures, and applications. Any progress in resolving this conjecture is  likely to reverberate through and impact a significantly larger swath of mathematics.
\pagebreak
\begin{center}
{\large REFERENCES}
\\
\end{center}
{\scriptsize \textit{Aand:} Stal Aanderaa,  Lars Kristiansen, Hans Kristian Ruud, Department of Mathematics, University of Oslo, A Preliminary Report on Search for Good Examples of Hall’s Conjecture, , January 20, 2014 (ArXiv:1401.4345v1)}
\\ \\
{\scriptsize \textit{Calv:} I. Jiminez Calvo, J. Herranz, G. Saez, A new algorithm to search for small nonzero $|x^{3} - y^{2}|$ values, Math. Comp. 76 (268) (2009) 2435-2444.}
\\ \\
\begin{scriptsize}
\textit{Danl:} Danilov, L.V.: The Diophantine equation $x^{3} - y^{2}$ and Hall's conjecture, Math. Notes Acad. Sci. USSR $\#32$ (1982), 617-618.
\end{scriptsize}
\\ \\
\begin{scriptsize}
\textit{Davn:} H. Davenport, On $f^{3}(t) - g^{2}(t)$, Norske Vid. Selsk. Forh. (Trondheim) 38 (1965), 86–87.
\end{scriptsize}
\\ \\
{\scriptsize \textit{Elks1:} Noam D. Elkies, Rational points near curves and small nonzero $|x^{3} - y^{2}|$ via lattice reduction. Algorithmic Number Theory. Proceedings of ANTS-IV; W. Bosma, ed.; Springer, 2000; pp. 33-63.}
\\ \\
{\scriptsize \textit{Elks2:} Noam D. Elkies, The ABC's of Number Theory, http://dash.harvard.edu/bitstream/handle/1/2793857/Elkies
\\ \\
{\scriptsize \textit{Gebl:} J. Gebel, A. Petho and H. G. Zimmer, On Mordell's equation, Compositio Math. 110 (1998), 335-367.}
\\ \\
{\scriptsize \textit{Gupt:} Gupta, Stebila, and Shantz, Integrating Elliptic Curve Cryptography into the Web’s Security Infrastructure}
\\ \\
{\scriptsize \textit{Malh:} K. Malhotra, S. Gardner, R. Patz, Implementation of Elliptic-Curve Cryptography on Mobile Healthcare Devices, Networking, Sensing and Control, April 2007 IEEE International Conference}
\\ \\
{\scriptsize \textit{Mazr1:} Mazur, Barry: Modular curves and the Eisenstein ideal, Inst. Hautes Etudes Sci. Publ. Math., No. 47, 1977, pp 33–186.}
\\ \\
{\scriptsize \textit{Mazr2:} Mazur, Barry: Questions about Number, http://www.math.harvard.edu/~mazur/papers/scanQuest.pdf}
\\ \\
{\scriptsize \textit{Mihl:} Preda Mihăilescu (2004). "Primary Cyclotomic Units and a Proof of Catalan's Conjecture". J. Reine angew. Math. 572 (572): 167–195. doi:10.1515/crll.2004.048. MR 2076124.}
\\ \\
{\scriptsize \textit{Mord:} L.J. Mordell, On the rational solutions of the indeterminate equations of the third and fourth degrees, Proc Cam. Phil. Soc. 21, (1922) p. 179.}
\\ \\
{\scriptsize \textit{Tidj:} Tijdeman, Robert (1976), "On the equation of Catalan", Acta Arithmetica 29 (2): 197–209, Zbl 0286.10013}
\\ \\
{\scriptsize \textit{Silv:} Silverman, Joseph, A quantitative version of Siegel's theorem: Integral points on elliptic curves and Catalan curves, J. Reine Angew. Math. 378 (1987), 60-100}
\\ \\
{\scriptsize \textit{Wald:} Michel Waldschmidt, Lecture on the abc conjecture and some of its consequences, Abdus Salam School of Mathematical Sciences (ASSMS), Lahore 6th World Conference on 21st Century Mathematics 2013}
\\ \\
{\scriptsize \textit{Weil:} A.	Weil, L'arithmétique sur les courbes algébriques, Acta Math 52, (1929) p. 281-315}
\\ \\ \\
{\scriptsize \textit{HTTP1:} https://oeis.org/A078933}
\\
{\scriptsize \textit{HTTP2:} $http://www.nsa.gov/business/programs/elliptic\_curve.shtml$}
\\
{\scriptsize \textit{HTTP3:} $http://blog.cloudflare.com/a-relatively-easy-to-understand-primer-on-elliptic-curve-cryptography$}
\\
{\scriptsize \textit{HTTP4:} $https://www.certicom.com/index.php/31-example-of-an-elliptic-curve-group-over-fp$}
\\
{\scriptsize \textit{HTTP5:} $http://www.networkworld.com/article/2224044/opensource-subnet/ecc-and-the-ca-security-council--making-ssl-and-the-web-safe-today-and-tomorrow.html$}
\\
{\scriptsize \textit{HTTP6:} http://oeis.org/A200216/internal}
\\
{\scriptsize \textit{HTTP7:} http://www.math.harvard.edu/~elkies/hall.html}
\\
{\scriptsize \textit{HTTP8:} http://www.ams.org/profession/prizes-awards/ams-supported/beal-prize}
\\
{\scriptsize \textit{HTTP9:} http://mathworld.wolfram.com/MordellCurve.html}
\\
{\scriptsize \textit{HTTP10:} http://oeis.org/A054504}
\\
{\scriptsize \textit{HTTP11:} http://web.math.pmf.unizg.hr/~duje/tors/rankhist.html}
\\
{\scriptsize \textit{HTTP12:} http://www.claymath.org/millennium-problems}
\end{document}